\documentclass[11pt]{article}

\usepackage{latexsym}
\usepackage{amsthm, amsmath, amssymb,latexsym,xypic}

\newtheorem{theorem}{Theorem}[subsection]

\newtheorem{lemma}[theorem]{Lemma}
\newtheorem{corollary}[theorem]{Corollary}
\newtheorem{proposition}[theorem]{Proposition}
\newtheorem{remark}[theorem]{Remark}
\newtheorem{definition}[theorem]{Definition}

\newtheorem*{teor31}{Theorem 3.1}

\newcommand{\nc}{\newcommand}
\nc{\bH}{{\mathbb H}}
\nc{\bA}{{\mathbb A}}
\nc{\bG}{{\mathbb G}}
\nc{\bC}{{\mathbb C}}
\nc{\bO}{{\mathbb O}}
\nc{\bI}{{\mathbb I}}
\nc{\bB}{{\mathbb B}}
\nc{\bY}{{\mathbb Y}}
\nc{\bK}{{\mathbb K}} 
\nc{\bX}{{\mathbb X}}
\nc{\bS}{{\mathbb S}}
\nc{\bE}{{\mathbb E}}
\nc{\bF}{{\mathbb F}}
\nc{\bZ}{{\mathbb Z}}
\nc{\bQ}{{\mathbb Q}}
\nc{\bN}{{\mathbb N}}
\nc{\bP}{{\mathbb P}}
\nc{\bL}{{\mathbb L}}
\nc{\bM}{{\mathbb M}}
\nc{\bT}{{\mathbb T}}
\nc{\bW}{{\mathbb W}}
\nc{\bU}{{\mathbb U}}
\nc{\bD}{{\mathbb D}}
\nc{\bJ}{{\mathbb J}}
\nc{\bV}{{\mathbb V}}
\nc{\bbZ}{{\mathbb Z}}
\nc{\bR}{{\mathbb R}}
\nc{\co}{{\nabla}}
\nc{\cu}{{\overline{\nabla}}}
\nc{\fr}{{\rightarrow}}
\newcommand{\la}{\longrightarrow}

\begin{document}

                                %
                                %
                                %
\title{Hermitian matrices and cohomology of K\"ahler varieties
} %
\author{Andrea Causin and Gian Pietro Pirola  \footnote{Partially supported by
1) PRIN 2005 {\em ``Spazi di moduli e
teorie di Lie"}; 
2) Indam (GNSAGA);
3) Far 2006 (PV):{\em ``Variet\`{a} algebriche, calcolo
algebrico, grafi orientati e topologici"}.}}
\date{}     

                                %
                                %
\maketitle
                                %
                                %

                                %
                                %


                                %
                                %
                                %
                                %
                                %
 \begin{abstract} {\em 
    We give some upper bounds on the dimension of the kernel of the cup product map 
 $H^{1}(X,\bC)\otimes H^{1}(X,\bC) \to H^{2}(X,\bC)$, where $X$ is a compact K\"ahler variety without Albanese fibrations.} \vskip3mm
\vskip 1mm
 {\setlength{\baselineskip}{0.8\baselineskip}
   \par}

\end{abstract}


                                %
                                %
                                %

 \section*{Introduction}

One of the special features of the K\"ahler geometry is the interplay
between topology and linear algebra. 
The goal of this paper is to give some applications of results in advanced linear algebra appearing in \cite{adams2} and \cite{frilib} to obtain bounds on the dimension of
the kernel of  the cup product mapping
 $$\phi \colon\bigwedge^2 H^{1}(X,\bC)  \to H^{2}(X,\bC),$$ in the
  case where $X$ is a compact K\"ahler variety admitting no Albanese fibration.
  
Recall that the notion of {\em no Albanese fibration, } for an $n-$dimensional 
variety $X$ as before, can be given by requiring that, for any $k<n$ and any independent 
$\beta_1,\dots ,\beta_k\in H^0(\Omega^1(X))$, the product $$\beta_1\wedge\dots\wedge\beta_k\in H^0(\Omega^k(X))=H^{k.0}(X)$$ is not zero. Also remark that 
$\wedge_{i=1}^k\beta_i=0$ gives rise to an integrable distribution and 
 to a foliation with closed leaves; this allows 
to define a fibration (see \cite{Catanese} or \ref{cdc} later in this paper), hence the
definition is consistent.
 
The main result we prove is:
\medskip

{\bf Theorem 1}
{\em Let $X$ be a compact K\"ahler 
variety without Albanese fibrations, and let
$\phi \colon\wedge^2H^{1}(X,\bC)\to H^{2}(X,\bC)$ be the cup product.
\begin{enumerate}
\item If  $q\leq 2n-1,$ then $\phi$ is injective.
\item If $ q=2n,$ then $\dim\ker\phi \leq  2c+3$ 
where $q=2^c(2b+1),$ and $b, c$  integers.
\item If $q=5$ and $n=2,$ then $\dim\ker\phi\leq  14.$ 
\end{enumerate} }
\medskip

The previous bounds are achieved thanks to the Hodge decomposition that we shortly recall. We have
$$H^{1}(X,\bC)=H^{1.0}\oplus H^{0.1};\ \ 
H^{2}(X,\bC)=H^{2.0}\oplus H^{1.1}\oplus H^{0.2},$$
hence $\phi$ defines maps:
\begin{enumerate}
\item $\phi^{2.0}: \bigwedge^2 
H^{1.0} \la H^{2.0}$ (and its dual $ \phi^{0.2}$);
\item $\phi^{1.1}: H^{1.0} \otimes H^{0.1} \la H^{1.1}. $
\end{enumerate}

The role of the assumption of no Albanese fibration becomes 
clear when dealing with the estimate of $\dim \ker \phi^{2.0}$
 (remark, for instance, that there are not decomposable elements in the kernel).
\medskip 

We look, then, for some upper-bound for  
  $\kappa=\dim\ker\phi^{1.1}.$
 Set   $H^{1.1}(X)_\bR= H^{1.1}(X)\cap H^2(X,\bR)$  and call $M\subset H^{1.0}\otimes H^{0.1}$ the subspace of forms invariant under complex conjugation.The space $M$ is naturally identified with the space $\bH_q$ of Hermitian $q\times q$ matrices and $\phi^{1.1}$ restricts to $\phi^{1.1}_\bR: M \rightarrow H^{1.1}(X)_\bR$; putting $K=\ker\phi^{1.1}_{\bR}$ we have  that $\kappa=\dim K.$
Assuming that $X$
 has no Albanese
fibration and using a positivity argument 
 we find restrictions on the
 signature and hence on the rank of the involved matrices:
if $A\in K,$ $A\neq 0,$ then the rank of $A$ must be $\geq 2n.$
In particular when $n=2q$ then $A$ is invertible and
from \cite{adams2} we achieve a very
 good bound on $\kappa.$ 
Write $q=2^c(2b+1),$ with $ b$ and $c$ integers; 
 we have $\kappa \leq 2c+1.$
The  basic remark is that the eigenvalues of
 $iA,$ with $i^2=-1,$ are not real, so if 
 $v\in S^{2q-1}\subset \bC^q$
 the tangential projection of $iAv$ on the unitary
  sphere defines a vector field
 not vanishing at any point.  
The result  is then based on  \cite{adams}.
\medskip

Unfortunately when $2n<q$ no good 
bound seems to be known in general.
However, at least for symmetric matrices,  
some very interesting works have been done; 
among the others, we suggest 
 \cite{ffl,frilib}. In a very strict sense this is a
{\em comeback of algebraic geometry}.
The basic idea is to consider the degeneracy locus 
$\{A: rank A\leq k\}$
of the matrices as real varieties,
then to study its intersection with real subspaces.
In particular in \cite{frilib}, 
after the natural projectivization, 
an elegant Lefschetz fixed point argument 
is used on a  suitable 
complex algebraic variety.  
This shows, in the case of $5\times 5 $ real symmetric 
matrices of rank lower than $4$, 
that the intersection is not empty.
  By performing the same kind of computations, for $n=2 $ and $q=5$ 
  we obtain 
  $\kappa\leq 8.$
\medskip

We would like to focus the importance of the mapping $\phi$ 
in the case of compact K\"ahler varieties.
On one side we have Castelnuovo-De Franchis-type theorem  
concerning Albanese fibrations  (see \cite{Catanese}).
On the other side, formality theorems \cite{deligneand} imply that 
the De Rham fundamental group $\pi_1(X)\otimes\bC$ is determined by $\phi$ \cite{wonderful}. 
We also recall works of Campana (see \cite{campana}) on K\"ahler nilpotent groups.
\medskip

The paper is organized in the following way. 
In section $1$ we recall the Albanese 
variety and the Castelnuovo-de Franchis
 theory as developed by  Catanese
\cite{Catanese}.
 In the second section we study the $1.1$
real forms and their connection with hermitian matrices.
This allows to prove the first two assertions 
of the above theorem. In section 3 we perform
 the computations necessary for the first case not covered 
in \cite{adams2}.\\

\noindent Both authors are delighted to thank Margherita 
for her helpful contribution and support.


\section{Variety of Albanese type}
\subsection{Albanese variety}
Let $X$ be a complex compact
 K\"ahler variety of dimension $n,$ and 
$$H^{p,q}=H^q(X,\Omega^p)$$ be the Hodge spaces. 
We have the Hodge decomposition:
$$H^n(X,\mathbb C)=
 \bigoplus_{p+q=n} H^{p,q}.$$
 Set $V=H^{1,0}$ and $H^{0,1}=\overline{V}$
   its conjugate.
Integration defines then  
$$j :  H_1(X,\mathbb Z) \la
V^{\ast},$$ where
$\ast$ stands for dual.
 Let $$Alb(X)= V^{\ast}/j(H_1(X,\mathbb Z))  $$ 
 be the Albanese variety. 
 The irregularity of $X$ is denoted by  $q_X=\dim V=\dim A.$ 
The choice of a base point $p\in X$ 
defines the Albanese map
$$ \alpha:X \la Alb(X).$$ 
\begin{definition} \label{type}
We say that  $X$ is of 
 Albanese type if $ \alpha$ is generically finite.
 We say that $X$ is of Albanese strict type if, 
 moreover, $\alpha$ is not surjective.
That is: $$\dim(\alpha(X))=\dim(X)<q_X.$$
\end{definition}
 
In the sequel we will assume $X$ 
of Albanese strict type. For our purposes,
thanks to the result of Campana,
we could also assume that $ \alpha$ 
 is generically one-to-one (see \cite{campana} and \cite[Ch.~2, Sect.~4]{wonderful}).


\subsection {Decomposable form and Albanese fibrations} \label{fibra}

We describe some  
Castelnuovo-de Franchis-type theorems.
With the previous notations, 
the cohomology map induced 
by the Albanese morphism:
$\alpha^{\ast}\colon H^k(A,\bC) \to H^k(X,\mathbb C)$
is a Hodge structure map. 
 We can moreover make the identifications: 
$H^k(A,\bC)\equiv \bigwedge^k H^1(X,\mathbb C),$
$V=H^{1.0}(A)\equiv H^{1.0}(X) $ 
and  
$H^{p.q}(A)=\bigwedge^pV\otimes\bigwedge^q\overline{V}.$
We have maps:
$\alpha^{p,q}\colon \bigwedge^pV\otimes \bigwedge^q\overline{V}\to H^{p.q}.$

\begin{definition} \label{Albanese-fibration}
Given $s>0,$ we say that a (rational map) $f:X\to Y$
 is an $s$-Albanese fibration if 
\begin{enumerate}
\item $Y$ is of Albanese strict type;
\item $\dim X-\dim Y=s.$
\end{enumerate}
When $ s=n-1,$ $ Y$ is a curve of genus 
 $g>1$ and $f$ is usually called an
  irregular pencil.
\end{definition}
\noindent We have (see \cite {Catanese}):
\begin{proposition}\label{cdc}
The following conditions are equivalent
\begin{enumerate}
\item $X$  has no $s$-Albanese fibration for  $s< n-k$
\item $\alpha^{k.0}$ is injective on the
decomposable forms: 
$$0\neq \beta_1\wedge \dots \wedge\beta_k\in H^{k,0}(X),$$
 $\beta_i$ independent.
\end{enumerate}
\end{proposition}

\begin{remark}
More precisely, Fabrizio Catanese in \cite{Catanese} 
 gave a one-to-one correspondence between 
 fibrations of Albanese type and maximal isotropic subspaces
  of the first cohomology group of $X.$
\end{remark}


\section{ Forms and matrices} 

\subsection{Real $1.1$ forms} \label{11}

As before, $X$
 is of strict Albanese type. 
We will consider  in  details the map:
\begin{equation} \label{a11}
\alpha^{1.1}\equiv \phi^{1.1}\colon H^{1.1}(A)
\equiv V \otimes \overline{V} \to H^{1.1}(X).
\end{equation}

\noindent We set 
\begin{equation}
\kappa=\dim \ker (\phi^{1.1}).
\end{equation} 
This is the $(1.1)$ part  of  $\phi\equiv \alpha^{2} : 
H^{2}(A,\bC)\equiv\bigwedge^2 H^{1}(X,\bC) \la H^{2}(X,\bC).$
Since $ \alpha^{2}$  is a piece of
a Hodge structure map,
the kernel of $ \alpha^{2}$ is defined over
 the rational numbers and a fortiori over the real numbers.
We have then a map:
$$ \alpha^{1.1}_{\bR}\colon H^{1.1}(A)_{\bR}
 \to H^{1.1}(X)_{\bR}\subset H^{2}(X,\bR)$$
and in particular $\kappa=\dim \ker (\alpha^{1.1}_{\bR}).$
 We can
identify $H^{1.1}(A)_{\bR}$  with the sesquilinear
forms on $V^{\ast}.$
 More explicitly, we fix a basis  $ \beta_j\  \ j=1,...,q$ of $V.$ 
An element of $H^{1.1}(A)_{\bR}$ has the form
$\Omega =i\sum_{j,s} a_{j,s} \beta_j\wedge \overline{\beta_s},$ $i^2=-1,$
 $a_{j,s}\in \bC:$  
 $\overline{\Omega}=\Omega,$
that is the matrix $A(\Omega)=(a_{j,s})$ is hermitian.
Let $\bH_q$ be the space of $q\times q$
 hermitian matrices.  
We may  define  
$\sigma\colon \bH_q\to H^{1.1}(X)_{\bR}$ as:  
\begin{equation} \label{rho}
\sigma(A)=i \sum_{j,s}a_{j,s}\ \beta_j\wedge\overline{\beta_s}, 
\end{equation} 
where $A=(a_{j,s}).$ 
We have $\kappa=\dim\ker\sigma.$
In particular, the rank and the signature of a
 form $\Omega\in H^{1.1}(A)_{\bR}$  are well defined.
\begin{definition} 
If $\Omega\in H^{1.1}(A)_{\bR}$ 
(respectively $A\in \bH_q$) has signature 
$(r,s)$ (and rank $r+s\leq q$),
we set $m(\Omega)= \min(r,s) $ (respectively $m(A)=\min(r,s)$).  
\end{definition}
The following proposition generalizes
 the elementary result explained
  in the introduction.
\begin{proposition} 
Let $k< n$ be an integer, and assume 
that $X$ has no $n-k$-Albanese fibration. 
Let $\Omega\in H^{1.1}(A)_{\bR}, \Omega\neq 0.$ 
If $m(\Omega)\leq k-1$ then $\alpha^{1.1}(\Omega)\neq 0.$
\end{proposition}
\begin{proof} 
Up to a change between 
$\Omega$ and $ -\Omega,$  
we may assume   $s=m(\Omega)$ where
$(r,s)$ is the signature of $\Omega,$ $s\leq k-1< n-1$. 
We may find a basis $\beta_i$ of $V$
such that
$$\Omega = i\sum_{j=1}^r   \beta_j\wedge\overline{\beta_j}- 
i\sum_{j=r+1}^{r+s}  \beta_j\wedge\overline{\beta_j}=\Omega^+ - \Omega^- .$$ 
Setting  $\varphi=\beta_{r+1}\wedge \dots\wedge \beta_{r+s}
 \in H^{s,0}(X) $ and 
$\Theta= \varphi\wedge\overline{\varphi},$ 
we compute:
$$  \Omega\wedge \Theta= 
i\sum_{j=1}^r  \beta_j\wedge\overline{\beta_j}\wedge \Theta=
(-1)^{s} i\sum_{j=1}^r  \beta_j \wedge\beta_{r+1}\wedge...\wedge
 \beta_{r+s}\wedge\overline{\beta_j \wedge\beta_{r+1}\wedge...\wedge \beta_{r+s}}.$$
Now, posing $ \varphi_j=  \beta_j \wedge\beta_{r+1}
\wedge...\wedge \beta_{r+s}\in H^{s+1,0}(X) $
and $\Theta_j= \varphi_j\wedge\overline{\varphi_j},$ we have:   
$$ \Omega\wedge \Theta=i (-1)^{s}\sum_j \Theta_j.$$
Assume by contradiction $\Omega\in \ker \alpha^{1.1}.$ 
It follows that $\Omega\wedge \Theta=0$
in $H^{s+1,s+1}(X).$
Fix  $\omega=i\sum_{j=1}^{q}\beta_j\wedge\overline{\beta_j};$ 
this is the pull-back of a K\"ahler form  on $A$ and
is positive on a Zariski open set of $X,$
since the Albanese map of $X$ 
is  generically finite. 
We get
$$0=\int_X \Omega\wedge \Theta\wedge \omega^{n-s-1}=
\sum_j \int_X i (-1)^{s} \Theta_j \wedge \omega^{n-s-1}.$$
All  terms have the same sign. 
It follows that 
$$ \Theta_j \wedge \omega^{n-s-1}=0;$$ 
this forces 
$\Theta_j=0$ and finally
$\varphi_j= \beta_j \wedge\beta_{r+1}\wedge...\wedge \beta_{r+s}=0.$
Since $s+1\leq k,$ we get a contradiction with Prop.~\ref{cdc}.
\end{proof}
We have then  the following:
\begin{corollary} \label{ranghi}
Assume that $X$ has no  Albanese fibration and
$\Omega\in \ker(\alpha^{1.1}) _{\bR},$
$\Omega \neq 0;$ then, $m(\Omega)\geq n-1$ and
  $ rank( \Omega) \geq 2n.$
\end{corollary}

\subsection{Hermitian matrices}

Let $\bH_q$ be the space of the $q\times q$ 
hermitian matrices.
Let $\bH_{q,m}\subset \bH_q$  be the subset 
of the matrices with rank bigger than $m-1:$
$$ \bH_{q,m}\in \{A\in \bH_{q}: rank(A)\geq m\}.$$

\begin{definition}
 Let $V\subset \bH_q$ be a real subspace. We
  say that $V$ has rank $\geq n$ ($rank(V)\geq n$) 
 if $V\setminus \{0\}\subset \bH_{q,n}.$
We set $d_{q,n} = \max_{rank (V)\geq n} \dim V.$
 \end{definition}
 \begin{proposition} \label{calcolo}
\begin{enumerate}
\item If  $q\leq m-1,$ then $ d_{q,m}=0.$ 
\item $d_{q,q}=2c+1,$ where $q=2^c(2b+1),$ with $b$ and $c$ integers.
\item $d_{5,4}\leq 8.$
\end{enumerate}
\end{proposition}
\begin{proof}
\begin{enumerate}
\item Obvious.
\item The  elements in $\bH_{q,q}$ 
are given by invertible hermitian matrices. Then, we come back  
to the hermitian case in \cite{adams2}.
This result was obtained as a consequence of \cite{adams}.
\item To be computed in the next section.
\end{enumerate}
\end{proof}

Recalling that 
$\kappa =\dim \ker \phi^{1.1}=\dim\ker \sigma$ (see \ref{rho}),  
we have the following:
\begin{proposition} \label{stima}
 Let $ n=\dim X,$ $ q=\dim H^{1.0}(X) $ 
 and $\kappa$ as before.  
Assume that $X$ has no Albanese 
fibration; then $\kappa \leq d_{q,2n}.$
\end{proposition}
\begin{proof} It follows from \ref{ranghi}. 
\end{proof}

Here we present some consequences of Prop.~\ref{calcolo}:
\begin{corollary} If  $q\leq 2n-1,$ then  
$\alpha^2\colon H^2(A,\bC) \to H^2(X,\bC)$ is injective.
\end{corollary}
\begin{proof}
Firstly one has that $\alpha^{1.1} \equiv\phi^{1.1}$
is injective. Then consider $\omega\in H^{2.0}(A)$ such that $\alpha^\ast (\omega)=0$; we show that $\omega=0$. Indeed, if it is not, we can find a basis $\beta_i$ of $V$ for which $\omega=\sum_{i=1}^k \beta_i\wedge\beta_{i+k}$, with $k<n$; taking the $k-1$ form $\phi=\wedge_{i=2}^k\beta_i$, we get $\alpha^\ast(\omega\wedge\phi)=0$ and consequently $\wedge_{i=1}^{k+1}\beta_i=0$ on $X$. By \ref{cdc}, this would give an Albanese fibration contradicting our assumptions.
\end{proof}

The previous result is essentially standard linear algebra. 
The first part of the following proposition is a consequence
of the hard topological result of Adams
\cite{adams,adams2}.  

\begin{proposition}\label{albadams}   
Assume than $X$ has no Albanese fibration and $q=2n.$   
Write $q=2^c(2b+1),$ $b$ and $c$ being integers.
Then we have
\begin{enumerate} 
\item $ \dim (\ker (\alpha^{1.1}))\leq 2c+1$ and
\item $\dim (\ker(\alpha^{2.0}))\leq  1.$
\end{enumerate}
The two inequalities above can be unified by saying that
 $\dim \ker (\alpha^{2})\leq 2c+3.$ Consequently it holds: $
b_2(X)\geq \dim Im \phi\geq q(2q-1) - 2c-3.$
\end{proposition}
\begin{proof}
\begin{enumerate}
\item  From \ref{stima} and \ref{calcolo}
we have  $\dim\ker \sigma\leq 2c+1. $
\item Arguing as in \ref{ranghi}, the nontrivial forms in
 $ \ker \alpha^{2.0}$ must be of maximal rank $n.$
Representing them as anti-symmetric matrices, the  
elements in $\ker  \alpha^{2.0}\setminus \{0\}$
are invertible. 
This is a complex space, and it follows that $\dim(\ker  \alpha^{2.0})\leq 1.$
\end{enumerate}
\end{proof}

\begin{proposition} Let $X$ be 
a compact algebraic variety 
without irregular pencils whose fundamental group admits a 
presentation with $\gamma$ generators and $\rho$ relations. 
 \begin{enumerate}
\item If  $q=2n,$  then  $\rho-\gamma\geq q(2q-3)-2c-3;$ moreover, 
if $X$ is a surface i.e. $n=2$, $c_2(X)\geq7.$   
\item If $q=5$ and $X$ is a surface, then $b_2(X) \geq 31,$ hence 
$\rho-\gamma\geq 31$ and $c_2(X)\geq 13.$ 
\end{enumerate}
\end{proposition}
\begin{proof}
Firstly recall that (see \cite[Th.~1.1]{amba} and \cite[Ch.~3]{wonderful}) 
the interplay 
between the cup product map and the fundamental group give rise to 
the following estimate: $\rho-\gamma\geq \dim Im (\phi)-2q.$ 
\begin{enumerate}
\item The estimate on $\rho - \gamma$ follows directly from the previous 
remark and \ref{albadams}. When $q=2n=4, $ we have $b_1(X)=b_3(X)=8$ 
and $b_2(X)\geq 21.$
\item To deduce $b_2(X)\geq 31$ one uses the third part of \ref{calcolo} 
together with the fact that 
$\dim Im\phi^{2.0} \geq 2q-3=7.$ The rest is as in the previous point.
\end{enumerate}
\end{proof}

\begin{remark}
In \cite{amba} the following estimate for a compact K\"ahler variety 
of any dimension with no irrational pencil is given: $\rho - \gamma \geq4q-7.$
\end{remark}

The following result was our historical motivation:
\begin{corollary} 
Let $X$ be a minimal projective
 algebraic surface with 
$q=4$ and $p_g=\dim H^{2.0}(X)=5.$ Let $K$ be the canonical
bundle of $X;$ then $ 16\leq  K^2\leq 17.$
\end{corollary}
\begin{proof} 
 In   \cite{bpn} it was proved that  $K^2\geq 16$
  and that if $X$ has an irregular pencil
then $K^2=16.$ When $X$ has no irregular pencil
Noether formula (see \cite{bpv,hirz})
$K^2+c_2(X) =12\chi_{\mbox{\small hol}}=
12(p_g-q+1)=24$ and
the inequality $c_2(X)\geq 7$ 
forces then $K^2\leq 17.$
 \end{proof}

\begin{remark}
In  the previous example the 
Miyaoka-Bogomolov-Yau inequality
(see for instance  \cite{bpv}) 
gives only $K^2\leq 9\chi_{\mbox{\small hol}}= 18.$ 
No examples of surfaces with 
$K^2=17$ are actually known.
\end{remark}

\section{Real  degeneracy loci.}
The present section is entirely dedicated to the proof of the third point in 
\ref{calcolo}. To this purpose, we adapt to hermitian matrices the work 
developed in \cite{frilib} in the case of real symmetric matrices, by thinking 
to the degeneracy locus $\{A \mbox{ hermitian, } rank A\leq m\}$ as a real variety.
The main theoric tool proved in \cite{frilib} is the following theorem, which 
is a consequence of the Hodge splitting and the Lefschetz fixed-point theorem: 

\begin{teor31}\label{libteo}
Let $V(\bR)\subset \bR\bP^n$ be an 
algebraic variety whose complexification
 $V\subset \bP^n $ is an irreducible variety 
 in codimension $m$. 
 Assume that the singular locus of $V$ 
has codimension at least $2r+1$ in $V$ 
and that, for a generic  $L\in Gr (m+2r+1, \bC^{n+1}),$  
the (topological) Euler characteristic $\chi(V\cap L)$
 is odd. Then, for any $P\in Gr (m+2r+1, \bR^{n+1})$, 
 $V(\bR)\cap P \neq \emptyset.$  
\end{teor31}
In our case, by choosing $V$ as the projectivization of the appropriate
 degeneracy locus (namely the $5\times 5$ complex matrices with rank 
 $3$ or less), and $L$ as a generic $8-$dimensional projective space,
  we show that the Euler characteristic of $V\cap L$ is odd,
   hence deducing that in the space of $5\times 5$ hermitian matrices 
   with rank not bigger than $3$ there are $8$ linearly independent elements. 
   This statement is equivalent to say that $d_{5,4}\leq 8.$ 
   About the calculation of $\chi(V\cap L)$, the arguments we 
   use are pretty standard ones and essentially concern Schubert 
   calculus of Grassmann spaces and Chern classes of projective bundles. 

\subsection{Chern classes of determinantal varieties.}
Endow the space of complex $q\times q$ matrices 
$M_q(\bC)$ with the real structure given by the
 (antiholomorphic) involution 
 $ A \mapsto \overline{A}^t$ 
 (where the bar stands for complex conjugation and $t$ for transposition); 
 with this choice, the real part of $M_q(\bC)$
  is the space $\bH_q$ of hermitian matrices. 
  Moreover, the real structure restricts to a real 
  structure over the irreducible affine variety 
  $V_{q,m}^a = \{A\in M_q(\bC) \ | \ rank(A)\leq m\}$
   and its real part $V_{q,m}^a(\bR)$ is exactly the
    complement in $\bH_q$ of the set $\bH_{q,m+1}$ 
    defined in the previous section. 
  
    Consider now the projectivizations 
$\bP^{q^2-1}= \bP(M_q(\bC))$ and
 $V_{q,m}=\bP(V_{q,m}^a),$ and set 
 $$ D_{q,m} = \min\{\dim P \ | \   V_{q,m}(\bR)\cap P \neq \emptyset \mbox{ for any linear }
  P \subset \bR\bP^{q^2-1} \}.$$ 
\begin{proposition} Let $d_{q,m}$ be
 as in the previous section. Then
  $d_{q,m+1}=D_{q,m}.$
\end{proposition}
\begin{proof} This should be clear since,
 by definition of $D_{q,m}$, there is 
 a linear space of projective dimension $D_{q,m} -1$
  entirely included in 
  $\bP(\bH_{q,m+1})=\bR\bP^{q^2-1}-V_{q,m}(\bR).$
\end{proof}
We would like, now, to make use of theorem \ref{libteo}.1:
we set $V=V_{5,3}\subset \bP^{24}$ 
with codimension $4$ and degree $50;$
its singular locus is $V_{5,2}$ with codimension 
$5$ in $V_{5,3}$ (see e.g. \cite{acgh}). 
Our goal is showing that $D_{5,3}\leq 8$, hence 
the only thing to prove is that, 
given a generic $L=\bP^8\subset \bP^{24}$, 
one has: 
\begin{equation}\label{caraeulero}\chi(V_{5,3}\cap \bP^8)\equiv 1 \mod (2).\end{equation}

\noindent To this purpose, consider the following
 resolution of $V_{5,3}$ with a smooth variety $Y$:
$$\xymatrix{  O(1) \ar@{}[r]|-{\supset}
 \ar[d]& O(1)|_{Y} \ar[d]  \ar@{}[r]|-{\simeq} & \mathcal O_\bP (1) \ar[d]  \\
		\bP(M_5(\bC))\times Gr(3,\bC^5) \ar@{}[r]|-{\supset}& Y=
		\{ ([A],W) \ | \ Im A\subset W \} \ar[dl]^{\pi_1} \ar[dr]_{\pi_2} \ar@{}[r]|-{\simeq} &
		 \bP(S^{\otimes 5}) \ar[d]  \\
		 V_{5,3} &  & Gr(3,\bC^5).} $$
 
 Here, $S\rightarrow Gr(3,\bC^5)$
 is the tautological bundle, $\bP(S^5)$ the
  projective bundle associated to 
  $S^{\otimes 5}$ and $\mathcal O_\bP(1)$
   and $O(1)$ are the duals to the tautological
    bundles over their respective projective bundles.
     Hyperplane sections of $O(1)|_Y$
      are exactly the same of those of $\mathcal O_\bP(1).$
      Remark that, outside $\pi_1^{-1}(V_{5,2})$, $Y$ is
       formed by couples of type $([A], ImA)$ hence  $\pi_1$ 
       is generically $1-1;$ remark also that the isomorphism 
       $Y\stackrel{\sim}{\to} \bP(S^5)$ is 
       given by the map $([A], W)\mapsto A\in
        Hom(\bC^5,W)\simeq W\otimes {\bC^\ast}^5 \simeq W^{\otimes 5}$.      
      
 For a vector bundle $F\rightarrow B$, let 
us denote its Chern classes by $ c_i(F)\in H^{2i} (B,\bZ) $ 
and its Chern polynomial $\sum_{k=0}^\infty c_k(F)t^k$ by $c_t(F).$

Computing $\chi(V_{5,3}\cap\bP^8)$ is 
the same as computing 
$\chi(Y\cap\bigcap_{i=1}^{16}H_i)$ 
with $H_i$ generic hyperplanes defined
 by sections of $O(1);$ since
 $Z=Y\cap\bigcap_{i=1}^{16}H_i$
  is a smooth complex $4$-dimensional 
 manifold, denoting by $T_Z$ 
 its tangent bundle and by $[Z]\in H_8(Z,\bZ)$
  its fundamental class, 
 we get  $$\chi(Z)= c_4(T_Z) [Z].$$ 

\noindent Following \cite{frilib} and \cite{hirz} we have:
\begin{proposition} Let $i\colon Z\to Y$ be the
 embedding and $h=c_1(O(1)|_Y)$;  then 
$$ c_t(T_Y) = i^\ast c_t(T_Y|_Z)(1+ht)^{-16}$$ and, 
denoting by $e_4$ the coefficient of $t^4$ in
 $c_t(T_Y)(1+ht)^{-16}:$
 \begin{equation}\label{caratteristica}
  \chi(Z)= h^{16}e_4 [Y].
 \end{equation}
\end{proposition}
\begin{proof} The first equality
 is a consequence of the exact 
 bundle sequence over $Z$: 
 $0\rightarrow T_Z \rightarrow i^\ast T_Y|_Z\rightarrow
  \bigoplus_{i=1}^{16} N_{H_i}|_Z \rightarrow 0$
   in which $N_{H_i}$ is the normal bundle of $H_i \subset Y$. 
   The second one follows as in (\cite[Ch.~9.2]{hirz}) since
    (the restriction of) $h$ represents a submanifold 
     in $Y\cap H_1\cap\dots\cap H_k$ for every $k.$   
\end{proof}
The diffeomorphism $Y\simeq \bP(S^5)$
 allows us to compute $c_t(T_Y)$ by means
  of the exact bundle sequences
   \begin{eqnarray}
    &0 \rightarrow  T_{\bP(S^5)/G}\rightarrow T_{\bP(S^5)} \rightarrow p^\ast T_G\rightarrow 0  
     \nonumber \\ & 0\rightarrow  \mathcal \bC \rightarrow S^5\otimes\mathcal O_\bP(1) \rightarrow Q\otimes \mathcal O_\bP(1)= T_{\bP(S^5)/G} \rightarrow 0 \nonumber 
   \end{eqnarray} 
   where the first one comes from the projection 
   $p:\bP(S^5)\rightarrow G=G(3,\bC^5)$ 
   and the second one is the tautological sequence 
   over $\bP(S^5 )$ tensorised with $\mathcal O_\bP(1).$ 
   Working out calculations we find (see \cite{frilib}):
\begin{equation}\label{prodottone} c_t(T_{\bP(S^5)})=
 c_t(T_G)c_t(T_{\bP(S^5)/G})=c_t(T_G)
 \left( \sum_{j=0}^{15} c_j(S^5) t^j(1+ht)^{15-j}\right).
 \end{equation}

\subsection{Computation of the Euler characteristic (modulo 2).}
We now dispose of  all the needed tools to complete
 the proof of the statement \ref{caraeulero}.\\
 Remark  that universal coefficient theorem implies  
that $H^\ast(B,\bZ)\otimes \bZ_2 = H^\ast(B, \bZ_2)$, 
where $B$ is either $G(3,\bC^5)$ or $\bP(S^5)$. 
Thanks to this, since we are only interested to the
 parity of our objects, we will perform any
  computation in $H^\ast(\ \cdot \ , \bZ_2); $
  moreover all polynomials in the Chern classes will be
  {1.1} truncated according to our real necessities.\\ 
Let us denote by $c_i$ the Chern classes of $S$.

\begin{proposition}\label{conto} 
Modulo $2$, we have:
\begin{enumerate} 
\item $c_t(S^5) = 1+c_1t +c_2t^2+c_3t^3 +c_1^4t^4+c_1^5t^5+c_1^4c_2t^6;$
\item $ c_t(T_G) = 1+c_1t+c_2t^2+c_3t^3 +c_2^2t^4+c_1c_2^2t^5+(c_1^4c_2+c_3^2)t^6;$
\item the  ring $H^\ast(\bP(S^5),\bZ_2)$ is $\bZ_2[c_1,c_2, h]$ 
together with the relations 
\begin{equation}\label{relazione} h^{15} = \sum_{j=1}^{15} c_j(S^5)h^{15-j} 
\qquad c_1c_2^2=0 \qquad c_1^4+c_1^2c_2+c_2^2=0
\end{equation}
\end{enumerate} 
\end{proposition}
\begin{proof} \begin{enumerate}
\item Follows directly from $c_t(S^5)=c_t(S)^5.$
\item Follows from the sequence 
$0\rightarrow S\otimes \tilde S\rightarrow \tilde S^5\rightarrow Q\otimes\tilde S 
=T_G  \rightarrow 0$ where $\tilde S$ is the dual of $S$ and 
$c_t(\tilde S)= c_{-t}(S)=c_t(S)$ since we 
are working in $\bZ_2.$ Hence: $c_t(T_G)=c_t(S\otimes\tilde S)c_t(S)^5.$
\item (For a more detailed treatment, see\cite{bottu}). 

The structure of the cohomology ring of 
a projective bundle $\bP(E)\rightarrow B$ 
is well known: 
$H^\ast(\bP(E))\simeq H^\ast(B)[h]/ (h^r+\sum_1^r c_i(E)h^{r-i})$ 
with $h$ the first Chern class of the
 canonical bundle over
 $\bP(E)$ and $r=rank(E)$.
 
 We have that $H^\ast(Gr(k,\bC^n),\bZ_2)=
 \bZ_2[c_1,\dots ,c_k]/(s_{n-k+1},\dots , s_n)$
  where $c_i$ are the Chern classes of the 
  tautological bundle $S\rightarrow Gr(k,\bC^n)$ 
  and $s_i$ those of the quotient bundle $\bC^n/S.$ 
  
In our context $s_3=c_1^3+c_3,$
$s_4=c_1^4+c_1^2c_2+c_2^2$ and $s_5=c_1^5+c_1^2c_3+c_1c_2^2;$
the statement follows from the elimination of $c_3$ 
and simplification of  
the expression of $s_5$.
\end{enumerate}
\end{proof}

\noindent Let  $e_4$ be the coefficient 
of $t^4$ in $c_t(T_{\bP(S^5)})(1+ht)^{-16}$;

\begin{lemma}\label{e_4} 
$e_4 =h^4+ c_1h^3+ c_1^2h^2+c_2h^2+c_1^3h+c_1c_2h +c_1^4$
 modulo $2.$
\end{lemma}
\begin{proof} From equation (\ref{prodottone}) we have
$$c_t(T_{\bP(S^5)})(1+ht)^{-16}=
c_t(T_G)(1+ht)^{-1}\left( \sum_{j=0}^{15} c_j(S^5) \left(\frac{t}{1+ht}\right)^j\right);$$ 
applying \ref{conto} and truncating polynomials to the $4$th 
degree, this becomes 
$$(1+c_1t+c_2t^2+c_1^3t^3 +c_2^2t^4)(1+ht+h^2t^2+h^3t^3+h^4t^4)
\sum_{j=0}^4c_j(S^5)(t+ht^2+h^2t^3+h^3t^4)^j.$$
The statement follows simply 
by working out calculations and by deleting couples of equal terms.
\end{proof}

\begin{proposition} \label{genera} 
The class $h^{16}e_4$ 
generates the top cohomology of $H^\ast(\bP(S^5),\bZ_2).$
\end{proposition}
\begin{proof} The proof of this statement is 
almost entirely based on the third point of \ref{conto}. 
The second and third relation of (\ref{relazione}) 
imply that $H^{12}(Gr(3,\bC^5),\bZ_2)$ is generated
 by any degree $6$ monomial in $c_1, c_2$, but $c_1^2c_2^2$. 
 
We claim that $H^{40}(\bP(S^5),\bZ_2)$
 is generated by $gh^{14}$ where $g$ is any 
  generator of $H^{12}(Gr(3,\bC^5),\bZ_2)$. 
  This can be seen directly by reducing the degree of 
  $h$ in any monomial of type $m(c_1,c_2)h^{14+k}$
   with the aid of the first relation of (\ref{relazione}). 
   As an example, consider $c_1^2c_2h^{16}:
   $ \begin{eqnarray}c_1^2c_2h^{16} & = c_1^2c_2h\ h^{15}
    =c_1^2c_2h  \  (c_1h^{14}+c_2h^{13} + \mbox{ higher order terms in  } c_1,c_2) =
    \nonumber \\ & = c_1^3c_2h^{15} + c_1^2c_2^2h^{14} =
      c_1^3c_2(c_1h^{14}+\dots) + 0 = c_1^4c_2 h^{14}.\nonumber\end{eqnarray}
Granting this, one finds (by analogous calculations)
 that the first and the three last terms in the expression 
 of $h^{16}e_4$ given by \ref{e_4} vanish,  
 while the three others do not. Hence,  $h^{16}e_4$ 
 generates $H^{40}(\bP(S^5),\bZ_2).$
\end{proof}

\begin{corollary} If $L$ is a generic $8$-dimensional
 linear subspace of $\bP^{24}$, the euler characteristic $\chi(L\cap V_{5,3})$ is odd.
\end{corollary}
\begin{proof} By equation (\ref{caratteristica}): 
$$\chi(L\cap V_{5,3})=\chi(Z)=h^{16}e_4[\bP(S^5)]
$$ and the latter term is not zero modulo $2$ thanks to \ref{genera} 
and Poincar\'e duality between $H^{40}(\bP(S^5),\bZ)$ and $H_{40}(\bP(S^5),\bZ).$ 
\end{proof} 

 \begin{remark} The following matrix, depending 
on $7$ real parameters, shows that $d_{5,4}=D_{5,3}\geq 7$, since $rank A\geq 4$ unless $\alpha=z=u=w=0.$
$$
A(\alpha,z,u,w)=\begin{pmatrix} \alpha 	& z 		& u  		& w 		& 0 \\
			      \bar z 	&\alpha	&\bar w	& -\bar u	& 0\\
			      \bar u	& w		& -\alpha	& z		& 0\\
			      \bar w	& -u		&\bar z	& -\alpha	& z \\
			      0		& 0		& 0		& \bar z 	& 0
\end{pmatrix} \qquad \alpha\in\bR, \ \  z,u,w \in \bC.
$$ 

To see this, let $A_k$ be the submatrix obtained from $A$  by elimination of the $k-$th row and column; $A_3$ is invertible unless $z=0$ or $|\alpha|=|z|$.
When $z=0$, $A_5^2= (\alpha^2+|u|^2 +|w|^2)I$ (see also \cite{adams2}); when $|\alpha|=|z|$, $\det A_1= |z|^2(|z|^2+|w|^2)$.
\end{remark}


                                %
                                %
                                %
                                %
                                %
                                %

                                %
                                %
                                %
                                %
                                %
                                %

\noindent
 {\sc Andrea Causin}\\
Dipartimento di Matematica, Universit\`a {\em``La Sapienza"} di Roma\\
P.le Aldo Moro 2, 00185 Roma, Italia\\
{\tt causin@mat.uniroma1.it}\\

\noindent {\sc Gian Pietro Pirola}\\
Dipartimento di Matematica, Universit\`a di Pavia\\
via Ferrata 1, 27100 Pavia, Italia\\
{\tt pirola@dimat.unipv.it}
\end{document}